\documentclass{article}

\usepackage{arxiv}

\usepackage[utf8]{inputenc} 
\usepackage[T1]{fontenc}    
\usepackage[numbers,sort&compress]{natbib} 
\usepackage{hyperref}       
\usepackage{url}            
\usepackage{booktabs}       
\usepackage{longtable}      
\usepackage{amsfonts}       
\usepackage{amsmath}
\usepackage{amssymb}
\usepackage{amsthm}
\usepackage{nicefrac}       
\usepackage{microtype}      
\usepackage{graphicx}
\usepackage{multirow}
\usepackage{tikz}
\usetikzlibrary{arrows.meta, positioning, fit, backgrounds}
\usepackage{listings}
\usepackage{algorithm}
\usepackage[noend]{algpseudocode}
\graphicspath{ {./} }

\newcommand{\dcut}{\delta}

\title{An open, reproducible branch-and-cut for the capacitated profitable tour
problem: a component study}

\author{
 Simon Spoorendonk \\
  Denmark \\
  \texttt{simon@spoorendonk.dk} \\
  \href{https://orcid.org/0009-0007-4304-6956}{ORCID:~0009-0007-4304-6956} \\
}

\begin{document}
\maketitle
\begin{abstract}
We present an open, reproducible branch-and-cut (B\&C) algorithm for the capacitated
profitable tour problem (CPTP) and its open $s$--$t$ path variant, the capacity-constrained
elementary shortest-path problem. The solver re-implements the formulation and cut families of
Jepsen et al.~\cite{jepsen2014cptp} on a fully open mixed-integer programming stack
(HiGHS~\cite{huangfu2018highs}), and adds bound-based preprocessing, domain propagation, and
reduced-cost variable fixing. We claim no new method; the contribution is twofold. First, an open,
reproducible artifact: to our knowledge the first branch-and-cut for this problem class on a fully
open stack, with the formulation, every separator, and all benchmark scripts released, so the
results below can be rerun and the solver reused and extended as a baseline. Second, a component
study on this common modern stack, benchmarked against a dynamic-programming/labelling reference,
that decomposes which components pay off and where the running time goes. We find that the
capacity-class cuts account for essentially the entire benefit (adding them to a connectivity-only
baseline lifts the number of instances solved from $52$ to $64$ of $76$ and shrinks the search tree
more than tenfold), while comb and rounded generalized-large-multistar cuts, reduced-cost fixing,
and bound-based propagation add nothing measurable. We also report a negative result: the
shortest-path-incompatibility (SPI) cut, a variant of the node-precedence inequalities of
García~\cite{garcia2009rcsp}, finds no violated inequality on any instance. The solver and all
experiments are released as open, reproducible software~\cite{cptp_software}.

\end{abstract}


\section{Introduction}
\label{sec:intro}

Let $G=(V,E)$ be an undirected graph with node set $V = N \cup \{0\}$, where $0$ is a designated
depot and $N$ is a set of customers. Each edge $e \in E$ has a cost $c_e$, and each node $i \in V$
has a profit $p_i \ge 0$ and a demand $d_i \ge 0$; the vehicle has capacity $Q$. The
\emph{capacitated profitable tour problem} (CPTP) asks for a simple cycle through the depot that
minimizes total edge cost minus collected node profit, subject to the total demand of the visited
nodes not exceeding $Q$. Because the objective trades travel cost against profit, optimal tours
visit a profitable subset of nodes rather than all of them, and the problem is naturally a
\emph{profitable tour} or \emph{prize-collecting} problem rather than a routing problem with a
fixed customer set.

CPTP is a special case of the elementary shortest-path problem with resource constraints (ESPPRC):
splitting the depot into a source and a sink and pushing node demands and profits onto arcs turns a
CPTP instance into an ESPPRC instance, though the reverse reduction does not hold in general
\cite{jepsen2014cptp}. ESPPRC and CPTP arise most prominently as the pricing subproblem in
column-generation algorithms for vehicle routing, where negative reduced costs on nodes create
exactly the profit-collection incentive modelled here. The problem is NP-hard, and when the
underlying cost structure admits negative-cost cycles (the regime that makes elementarity
binding), it is the hard core of routing pricing.

Two lines of attack dominate the literature. \emph{Dynamic-programming / labelling} algorithms
extend partial paths under dominance and resource bounds; they are the workhorses inside modern
column generation and are extremely effective when the resource (here, capacity) keeps the number
of non-dominated labels small \cite{feillet2004espprc, lozano2016espprc, puglieseguerriero2012rcespp}. Their
performance degrades as elementarity, rather than the resource, becomes the binding difficulty.
\emph{Branch-and-cut} (B\&C) algorithms instead solve an integer-programming formulation and
separate valid inequalities. For the capacitated, negative-cycle regime, the closest prior art is
the doctoral work of García~\cite{garcia2009rcsp}, who developed a B\&C for resource-constrained
shortest paths with shortest-path-bound cuts (node-precedence and conflict-graph inequalities),
preprocessing, and branching, and the CPTP-specific B\&C of Jepsen et al.~\cite{jepsen2014cptp},
which catalogues and separates a large family of valid inequalities and shows B\&C to be competitive
with and complementary to labelling on capacitated instances. Related cutting-plane work includes
Avella et al.~\cite{avella2004rcsp} and Horváth and Kis~\cite{horvath2016rcspp}; for the pure
elementary shortest-path problem without a resource, B\&C formulations were studied by
Taccari~\cite{taccari2016espp} and Drexl and Irnich~\cite{drexl2014espp}, though these exercise
only connectivity cuts and not the capacity machinery. Subsequent exact work has extended rather
than superseded this line (Lera-Romero and Miranda-Bront~\cite{leraromero2021tdcptp} give a
branch-and-cut for a time-dependent, resource-constrained profitable tour), but we are not aware of
a later exact method that improves on Jepsen et al.~\cite{jepsen2014cptp} for the base CPTP. The
ESPPRC pricing problem itself has meanwhile advanced through bidirectional~\cite{righini2006symmetry}
and bucket-graph labelling~\cite{sadykov2021bucket} embedded in branch-price-and-cut
frameworks~\cite{costa2019bpc, pessoa2020vrpsolver}; our solver targets the standalone subproblem
rather than a full routing solver.

The methodological space is therefore well charted: that branch-and-cut can match dynamic
programming on negative-cycle capacitated paths was established for the CPTP, with capacity-class
cuts, by Jepsen et al.~\cite{jepsen2014cptp}, and for the broader resource-constrained shortest
path, with shortest-path-bound cuts and fixing, by García~\cite{garcia2009rcsp}. We do not claim a
new algorithm. Instead, the contribution of this paper is empirical and infrastructural:

\begin{itemize}
\item \textbf{An open, reproducible artifact.} To our knowledge this is the first branch-and-cut
  for this problem class built end-to-end on a fully open MIP stack (HiGHS~\cite{huangfu2018highs},
  C++23), with the formulation, all separators, preprocessing, propagation, and heuristics released
  as open software together with every benchmark script. It is a rerunnable, reusable baseline that
  others can build on and extend, for example to the per-arc, directed models discussed in
  Section~\ref{sec:instances}. The solver also handles the open capacity-constrained elementary
  shortest-path ($s$--$t$) variant.
\item \textbf{A component study.} We report a systematic ablation of the cut families, of domain
  propagation, and of cut-based variable fixing on a common modern stack: a decomposition of which
  components pay off and where the running time goes that, to our knowledge, has not previously been
  reported for this problem.
\item \textbf{A negative result.} The shortest-path-incompatibility (SPI) cut implemented in
  the solver (a variant of the node-precedence and conflict-graph inequalities of
  García~\cite{garcia2009rcsp}, with roots in the shortest-path-bound preprocessing of Aneja et
  al.~\cite{aneja1983shortest}) yields no measurable benefit on our instances. We document this so
  that others are not drawn to a dead cut.
\end{itemize}

The remainder of the paper is organized as follows. Section~\ref{sec:model} states the integer
programming model and the $s$--$t$ path variant. Section~\ref{sec:bac} describes the branch-and-cut
algorithm: the cut families and their separation, and the preprocessing, propagation, and fixing
that surround it. Section~\ref{sec:study} reports the computational study: the comparison with
dynamic programming, the gating ablation, and the SPI negative result. Section~\ref{sec:conclusion}
concludes.

\section{Problem definition and model}
\label{sec:model}

We use the following shorthand, following Jepsen et al.~\cite{jepsen2014cptp}. For $S \subseteq V$,
$\dcut(S)$ denotes the set of edges with exactly one endpoint in $S$ (the cut boundary), and $E(S)$
the set of edges with both endpoints in $S$. Singletons are abbreviated, so $\dcut(i)$ means
$\dcut(\{i\})$. The binary variable $x_e$ indicates whether edge $e$ is used and $y_i$ whether node
$i$ is visited; a fractional LP solution is written $(\bar{x}, \bar{y})$. One-customer tours (the
depot and a single customer) can be enumerated a priori in $O(|N|)$ and supplied as upper bounds;
the formulation still represents them directly, by letting a depot-incident edge take value $2$.

The tour variant fixes the depot in the solution ($y_0 = 1$) and allows the depot-incident edges to
take the value $2$ noted above. The integer programming model is

\begin{align}
\min \quad & \sum_{e \in E} c_e x_e - \sum_{i \in V} p_i y_i \label{eq:obj}\\
\text{s.t.}\quad
& \sum_{e \in \dcut(i)} x_e = 2 y_i && \forall i \in V \label{eq:degree}\\
& y_0 = 1 \label{eq:depot}\\
& \sum_{i \in V} d_i y_i \le Q \label{eq:cap}\\
& \sum_{e \in \dcut(S)} x_e \ge 2 y_i && \forall S \subset V \setminus \{0\},\ \forall i \in S \label{eq:gsec}\\
& x_e \in \{0,1\} && \forall e \in E \setminus \dcut(0) \label{eq:xbin}\\
& x_e \in \{0,1,2\} && \forall e \in \dcut(0) \label{eq:xdepot}\\
& y_i \in \{0,1\} && \forall i \in V. \label{eq:ybin}
\end{align}

The objective~\eqref{eq:obj} minimizes travel cost minus collected profit. The degree
constraints~\eqref{eq:degree} couple edge and node variables: a visited node has degree two, an
unvisited node degree zero. Constraint~\eqref{eq:depot} forces the depot into the tour, and the
knapsack constraint~\eqref{eq:cap} enforces the capacity. The generalized subtour elimination
constraints (GSEC)~\eqref{eq:gsec} guarantee connectivity: any visited node $i \in S$ must be linked
to the rest of the graph across $\dcut(S)$. There are exponentially many GSEC, so they are separated
dynamically rather than enumerated (Section~\ref{sec:bac}).

Following Jepsen et al.~\cite{jepsen2014cptp}, we name the polytopes that organize the valid
inequalities. The \emph{circuit polytope} $P_C$ is defined by the degree
constraints~\eqref{eq:degree}, the GSEC~\eqref{eq:gsec}, and the integrality~\eqref{eq:xbin}--\eqref{eq:ybin};
the \emph{knapsack polytope} $P_K$ by the capacity constraint~\eqref{eq:cap} and
integrality~\eqref{eq:ybin}; and the \emph{node-capacitated circuit polytope}
$P_N = P_C \cap P_K$ is the feasible region of the CPTP. Cut families are catalogued by which
polytope justifies their validity.

\paragraph{The open $s$--$t$ path variant.}
When a distinct source $s$ and target $t$ are specified, the model solves the open
capacity-constrained elementary shortest path from $s$ to $t$: constraint~\eqref{eq:depot} and the
two-valued depot edges are dropped, the degree of $s$ and $t$ is fixed to one, and the GSEC are
imposed relative to the terminals. This variant is the natural comparator for ESPPRC labelling
algorithms and is used as such in Section~\ref{sec:study}.

\section{Branch-and-cut}
\label{sec:bac}

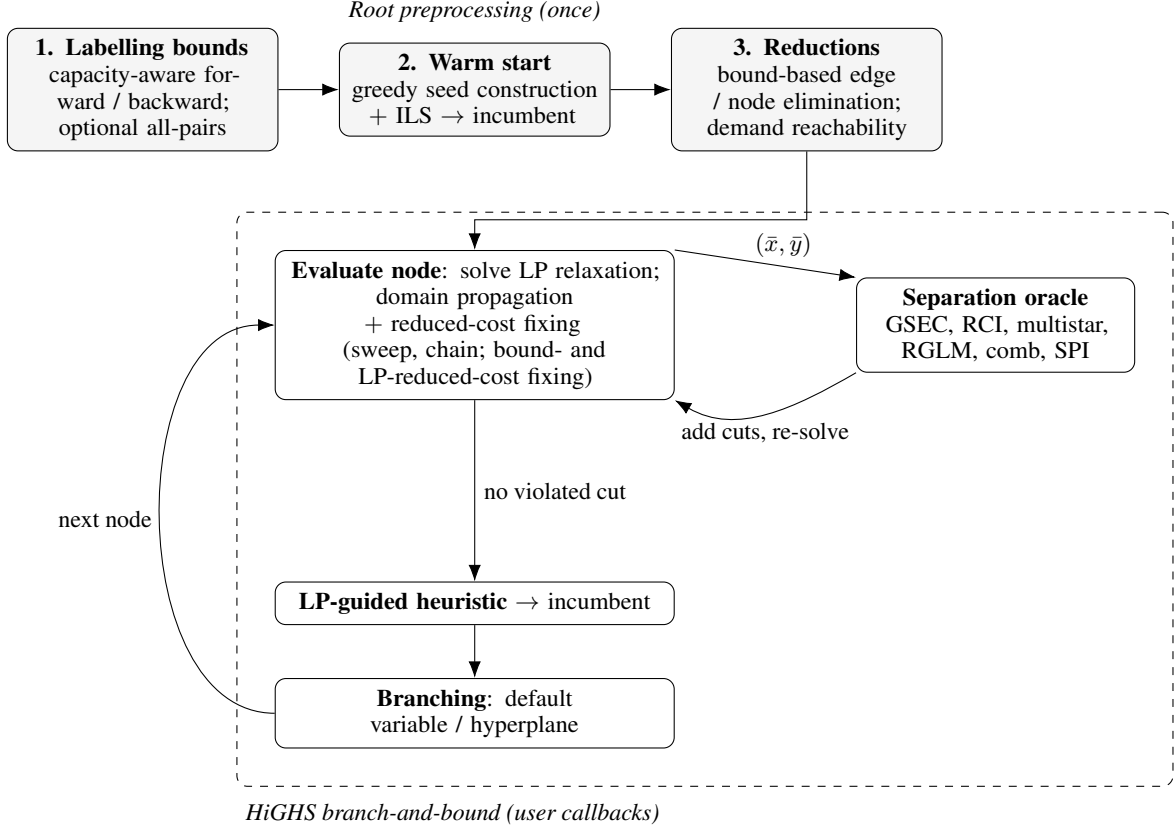
\begin{figure}[htbp]
\centering
\resizebox{\ifdim\width>\linewidth\linewidth\else\width\fi}{!}{%
\begin{tikzpicture}[
  font=\small,
  >={Latex[length=2.2mm]},
  block/.style={draw, rounded corners, align=center, inner sep=4pt},
  pre/.style={block, fill=black!4, text width=33mm, font=\footnotesize},
  node distance=7mm,
]
\node[pre] (p1)
  {\textbf{1. Labelling bounds}\\ capacity-aware forward / backward; optional all-pairs};
\node[pre, right=8mm of p1] (p2)
  {\textbf{2. Warm start}\\ greedy seed construction $+$ ILS $\rightarrow$ incumbent};
\node[pre, right=8mm of p2] (p3)
  {\textbf{3. Reductions}\\ bound-based edge / node elimination; demand reachability};
\draw[->] (p1) -- (p2);
\draw[->] (p2) -- (p3);
\node[font=\footnotesize\itshape, above=1.5mm of p2] {Root preprocessing (once)};

\node[block, below=15mm of p2, text width=50mm] (lp)
  {\textbf{Evaluate node}: solve LP relaxation;\\ domain propagation $+$ reduced-cost fixing\\ {\footnotesize(sweep, chain; bound- and LP-reduced-cost fixing)}};
\node[block, below=24mm of lp, text width=50mm] (heu)
  {\textbf{LP-guided heuristic} $\rightarrow$ incumbent};
\node[block, below=of heu, text width=50mm] (br)
  {\textbf{Branching}: default variable / hyperplane};
\node[block, right=24mm of lp, text width=34mm] (sep)
  {\textbf{Separation oracle}\\ GSEC, RCI, multistar,\\ RGLM, comb, SPI};

\draw[->] (p3.south) |- ([yshift=4mm]lp.north) -- (lp.north);
\draw[->] (lp.north east) -- node[above, pos=0.6, font=\footnotesize]{$(\bar x,\bar y)$} (sep.north west);
\draw[->] (sep.south west) to[out=210, in=-30]
  node[below, pos=0.5, font=\footnotesize]{add cuts, re-solve} (lp.south east);
\draw[->] (lp) -- node[right, pos=0.5, font=\footnotesize]{no violated cut} (heu);
\draw[->] (heu) -- (br);
\draw[->] (br.west) to[out=180, in=180] node[left, pos=0.5, font=\footnotesize]{next node} (lp.west);

\begin{scope}[on background layer]
\node[draw, dashed, rounded corners, inner sep=5mm,
      fit=(lp)(sep)(heu)(br)] (bbox) {};
\end{scope}
\node[anchor=north west, font=\small\itshape] at ([yshift=-1mm]bbox.south west)
  {HiGHS branch-and-bound (user callbacks)};
\end{tikzpicture}%
}
\caption{Where each component acts. Root preprocessing runs once (top); the separation oracle,
domain propagation with reduced-cost fixing, the LP-guided heuristic, and branching are HiGHS
callbacks invoked per node. Domain propagation and reduced-cost fixing run \emph{with} each LP
solve (and may tighten bounds and trigger a re-solve); the separation oracle then adds the most
violated cuts and the LP is re-solved, looping until no cut is found, after which the node yields to
the heuristic and branching. The cut families, the propagation and reduced-cost fixing, and the
all-pairs/SPI option are exactly the components toggled in the ablation of
Section~\ref{sec:ablation}.}
\label{fig:flow}
\end{figure}

The model~\eqref{eq:obj}--\eqref{eq:ybin} is solved with the HiGHS mixed-integer
solver~\cite{huangfu2018highs}, extended through callbacks for user cut separation, domain
propagation, and branching. Presolve is disabled to keep a stable column mapping for the callbacks.
We first describe the cut families and their separation, then the preprocessing, propagation, and
fixing that surround the search. Figure~\ref{fig:flow} shows where each component acts.

\subsection{Valid inequalities and their separation}
\label{sec:cuts}

All separators operate on the fractional \emph{support graph} $G_{\bar x}$ of the current LP
solution $(\bar{x}, \bar{y})$: each edge $e$ with $\bar{x}_e > \varepsilon$ ($\varepsilon =
10^{-6}$) contributes both arcs $(u,v),(v,u)$ of capacity $\bar{x}_e$, so that the value of any
node-cut equals $\bar{x}(\dcut(S))$. A Gomory--Hu cut tree~\cite{gusfield1990gomoryhu}, built by
Gusfield's construction in $|V|-1$ maximum-flow computations~\cite{dinic1970maxflow}, is computed
once per separation round, rooted at the depot, and shared by the connectivity- and capacity-class
families; each tree edge encodes a minimum $i$--depot cut and the cuts \emph{nested} along a
root-path are read off at no extra flow cost. The complete per-round procedure is
Algorithm~\ref{alg:sep}. Throughout, for $e \in \dcut(S)$ we write $o(e)$ for the endpoint of $e$
\emph{outside} $S$, and $d(S) = \sum_{i \in S} d_i$. A separator emits a cut only when its violation
exceeds a tolerance (default $0.1$); at most $K = 10$ cuts, those with the largest violation, are
added per round across all families.

\paragraph{Generalized subtour elimination (GSEC).}
For each customer $i$ with $\bar{y}_i > 0$ the separator reads the minimum $i$--depot cut $(f_i, S_i)$
from the tree ($S_i$ the shore containing $i$, excluding the depot) and the inequality~\eqref{eq:gsec}
is violated when
\begin{equation}
\bar{x}(\dcut(S_i)) \;=\; f_i \;<\; \kappa\,\bar{y}_i,
\qquad \kappa = \begin{cases} 2 & \text{tour, or } t \notin S_i,\\ 1 & \text{path with } t \in S_i, \end{cases}
\label{eq:gsec-sep}
\end{equation}
the value $\kappa=1$ holding in the path variant when $S_i$ contains the terminal $t$, since the path
need only cross $\dcut(S_i)$ once. Substituting the degree equations $\bar{x}(\dcut(j))=2\bar{y}_j$
yields the algebraically equivalent \emph{inside} form $x(E(S)) \le \sum_{j \in S} y_j - y_i$, valid
only for terminal-free $S$; the solver emits whichever of the two forms is sparser. A strong-component
separation with a better worst case is given by Drexl~\cite{drexl2013sec}.

\paragraph{Rounded capacity inequalities (RCI).}
Capacity implies that any set $S$ requires at least $\lceil d(S)/Q \rceil$ ``loads'' crossing its
boundary, the classic rounded-capacity cut $\sum_{e \in \dcut(S)} x_e \ge 2\lceil d(S)/Q \rceil$. We
separate the residual-demand strengthening of Jepsen et al.~\cite{jepsen2014cptp}: with $Q_r = d(S)
\bmod Q$ and $k = \lceil d(S)/Q \rceil$, and only when $Q_r > 0$ and $k \ge 2$,
\begin{equation}
\sum_{e \in \dcut(S)} x_e \;\ge\; \frac{2}{Q_r}\sum_{i \in S} d_i\,y_i \;+\; 2\!\left(k - \frac{d(S)}{Q_r}\right).
\label{eq:rci}
\end{equation}
The right-hand side equals $2k$ when all of $S$ is selected ($y_i \equiv 1$) and relaxes
monotonically with the served demand $\sum_{i\in S} d_i y_i$, so~\eqref{eq:rci} coincides with the
classic bound at full selection and is its valid node-coupled generalization (the classic constant
$2\lceil d(S)/Q\rceil$ need not be valid once nodes of $S$ may go unvisited). Candidate sets are the nested cuts on each customer's tree
root-path, each refined by an alternating add/drop hill-climb that grows or shrinks $S$ across
boundary edges to maximize violation; up to $\max(5, |V|/5)$ RCI are kept per round.

\paragraph{Multistar and generalized large multistar (GLM).}
The multistar (capacity) inequality
$\sum_{e \in \dcut(S)} x_e \ge \frac{2}{Q}\sum_{i \in S} d_i y_i$
strengthens the capacity bound by placing \emph{variables}, not constants, on the right-hand side.
The generalized large multistar of Letchford and Salazar-González~\cite{letchford2006projection}
(in the one-vehicle CPTP form of Jepsen et al.~\cite{jepsen2014cptp}) additionally credits, on the
boundary, the demand of the outside endpoint of each crossing edge:
\begin{equation}
\sum_{e \in \dcut(S)} \!\left(1 - \frac{2\,d_{o(e)}}{Q}\right) x_e \;\ge\; \frac{2}{Q}\sum_{i \in S} d_i\,y_i .
\label{eq:glm}
\end{equation}
For each target customer the single minimum cut $S$ from the tree is tested against~\eqref{eq:glm}.

\paragraph{Rounded GLM (RGLM).}
The rounded generalized large multistar (eq.~(24) of Jepsen et al.~\cite{jepsen2014cptp}) ceiling-rounds
the GLM. Writing $\hat{d} = \sum_{i} d_i$ for the total eligible demand, $\alpha = 2\hat{d} - d(S)$,
$r = \alpha \bmod Q$, and $k = \lceil \alpha/Q \rceil$ (separated when $r > 0$, $k \ge 2$),
\begin{equation}
\sum_{e \in \dcut(S)} \!\left(1 - \frac{2\,d_{o(e)}}{r}\right) x_e \;\ge\; \frac{2}{r}\sum_{i \in S} d_i\,y_i \;+\; 2k - \frac{2\alpha}{r}.
\label{eq:rglm}
\end{equation}
The GLM's domination of the plain capacity inequality does not survive rounding, so both are kept.
RGLM is built on the same minimum-cut sets as the GLM/capacity separators.

\paragraph{Comb inequalities.}
A comb has a handle $H$ and an odd number $T \ge 3$ of teeth; here each tooth is a single boundary
edge $e_t$ with inside endpoint in $H$ and outside endpoint $o_t \notin H$. In the node-coupled CPTP
form, adapting Bauer's circuit-polytope combs~\cite{bauer1997circuit} as in Jepsen et
al.~\cite{jepsen2014cptp}, the inequality is
\begin{equation}
\sum_{e \in E(H)} x_e \;+\; \sum_{t=1}^{T} x_{e_t} \;-\; \sum_{j \in H} y_j \;\le\; \frac{T-1}{2}.
\label{eq:comb}
\end{equation}
Separation is heuristic: a handle is grown by breadth-first search from the depot over edges with
$\bar{x}_e$ above a threshold (tried at $0.5$ and $0.3$), candidate teeth are the boundary edges
ranked by $\bar{x}_e - \bar{y}_{o}$, a maximal set of node-disjoint teeth is taken greedily, and the
count is reduced to the nearest odd number. In the $s$--$t$ variant the two terminals are merged into
a single handle root. At most two combs (one per threshold) are produced per round.

\paragraph{Shortest-path-incompatibility (SPI) inequalities.}
The SPI separator derives node-incompatibility cuts from all-pairs bounds together with the incumbent.
Let $L(S)$ be a lower bound on the cost of any feasible tour (or $s$--$t$ path) that visits all of
$S$, computed by an exact Held--Karp dynamic program~\cite{heldkarp1962dp} over $S$ between the
terminals ($O(2^{|S|}|S|^2)$, capped at $|S| \le 15$). If $L(S)$ exceeds the incumbent upper bound
$\overline{z}$, the nodes of $S$ cannot all be selected, giving the cover cut
\begin{equation}
\sum_{i \in S} y_i \;\le\; |S| - 1
\qquad\text{whenever}\qquad L(S) > \overline{z},
\label{eq:spi}
\end{equation}
with the pair case $y_i + y_j \le 1$. The separator seeds from violated pairs, grows each seed by
adding high-$\bar{y}$ nodes while $L(S) > \overline{z}$, shrinks to a minimal infeasible set, and
lifts further nodes whose every one-for-one swap into $S$ stays infeasible (so the right-hand side is
unchanged).

This family is \emph{not} a new class of inequality. It is a variant of the node-precedence and
conflict-graph inequalities of García~\cite{garcia2009rcsp}, whose validity rests on the
shortest-path-bound feasibility arguments that descend from the preprocessing of Aneja et
al.~\cite{aneja1983shortest}. The only element specific to the present solver is the use of a
Held--Karp bound in place of a plain shortest-path bound: a tightening of a known cut, not a new
inequality. Section~\ref{sec:study} reports that it gives no empirical benefit.

\begin{algorithm}[htbp]
\caption{Separation round at an LP solution $(\bar x,\bar y)$ with incumbent $\overline{z}$}
\label{alg:sep}
\begin{algorithmic}[1]
\State build support graph $G_{\bar x}$: arcs $(u,v),(v,u)$ of capacity $\bar x_e$ for every $e$ with $\bar x_e > \varepsilon$
\State build Gomory--Hu cut tree $\mathcal{T}$ on $G_{\bar x}$, rooted at the depot \Comment{$|V|-1$ max-flows}
\State $\mathcal{C} \gets \emptyset$ \Comment{cut pool}
\For{each customer $i$ with $\bar y_i > 0$} \Comment{run in parallel over families}
  \State $(f_i, S_i) \gets$ minimum $i$--depot cut read from $\mathcal{T}$
  \If{$f_i < \kappa \bar y_i$} add GSEC~\eqref{eq:gsec-sep} for $S_i$ to $\mathcal{C}$ \EndIf
  \State test multistar/GLM~\eqref{eq:glm} and RGLM~\eqref{eq:rglm} on $S_i$; add if violated
  \For{each nested set $S'$ on the tree path $i \to$ depot}
    \State refine $S'$ by add/drop hill-climb; if RCI~\eqref{eq:rci} violated, add to $\mathcal{C}$
  \EndFor
\EndFor
\State grow comb handles from the depot at thresholds $\{0.5, 0.3\}$; add violated combs~\eqref{eq:comb}
\If{all-pairs bounds available} run SPI~\eqref{eq:spi}; add violated cover/pair cuts to $\mathcal{C}$ \EndIf
\State sort $\mathcal{C}$ by violation; \Return the $K$ most violated cuts
\end{algorithmic}
\end{algorithm}

\subsection{Preprocessing, propagation, and fixing}
\label{sec:fixing}

Around the cut loop, the solver applies bound-based reductions that descend from the
shortest-path-bound preprocessing of Aneja et al.~\cite{aneja1983shortest} and
García~\cite{garcia2009rcsp}; we treat them as implementation, not new ideas.

\paragraph{Capacity-aware labelling bounds.}
A label at node $v$ carries two resources, accumulated net cost (edge cost minus collected profit)
and accumulated demand. Starting from a root, labels are extended along edges, rejecting any
extension whose demand would exceed $Q$, and pruned by Pareto dominance on (cost, demand)
(Algorithm~\ref{alg:label}). The pass returns, for every node $v$, the minimum net cost $f_v$ of a
capacity-feasible walk from the source to $v$; a backward pass from the target gives $b_v$ (for tours
the two coincide). The labelling is \emph{non-elementary} (it tracks no visited set, only the two
resources, so a node may repeat), hence $f_v, b_v$ are valid \emph{lower bounds}, possibly loose;
this relaxation is what makes them cheap, and bounds the label set so the pass runs under a fixed
pop budget. An optional all-pairs variant runs the forward labelling from every node, yielding a
matrix of pairwise bounds $d(\cdot,\cdot)$ that the chain propagation below uses to route through any
pair of edges.

\begin{algorithm}[htbp]
\caption{Capacity-aware labelling bound from a root $r$ (forward pass; backward is symmetric from the target)}
\label{alg:label}
\begin{algorithmic}[1]
\Require graph $G$, edge costs $c$, demands $d$, profits $p$, capacity $Q$, pop budget $B$
\Ensure $f_v$, a lower bound on the net cost of a capacity-feasible walk $r \to v$, for all $v$
\State $f_v \gets +\infty$, $\mathcal{L}_v \gets \emptyset$ for all $v$ \Comment{$\mathcal{L}_v$: Pareto label set at $v$}
\State $\ell_0 \gets (\,c{:}\,{-}p_r,\ q{:}\,d_r\,)$;\quad $f_r \gets -p_r$;\quad push $(r,\ell_0)$ to queue; $\mathcal{L}_r \gets \{\ell_0\}$
\State $\mathit{pops} \gets 0$
\While{queue $\neq \emptyset$ \textbf{and} $\mathit{pops} < B$}
  \State $(u,(c_u,q_u)) \gets$ pop; $\mathit{pops} \gets \mathit{pops}+1$
  \For{each edge $e=(u,v)$}
    \State $q' \gets q_u + d_v$; \textbf{if} $q' > Q$ \textbf{then continue} \Comment{capacity}
    \State $c' \gets c_u + c_e - p_v$
    \If{$\exists\,(c'',q'') \in \mathcal{L}_v$ with $c'' \le c'$ and $q'' \le q'$} \textbf{continue} \Comment{dominated}
    \EndIf
    \State remove from $\mathcal{L}_v$ every label dominated by $(c',q')$; add $(c',q')$ to $\mathcal{L}_v$
    \State $f_v \gets \min(f_v, c')$; push $(v,(c',q'))$ to queue
  \EndFor
\EndWhile
\State \textbf{if} $\mathit{pops} = B$ \textbf{then} \Return \textsc{Abort} \Comment{budget hit: disable bound-based features}
\State \Return $f$
\end{algorithmic}
\end{algorithm}

\paragraph{Bound-based elimination.}
Let $\overline{z}$ be the incumbent upper bound and $\kappa = p_0$ for a tour (the depot profit,
which both passes net out) or $\kappa = 0$ for a path. Any capacity-feasible solution using edge
$e=(u,v)$ costs at least
\begin{equation}
\ell(e) = \min\bigl(f_u + c_e + b_v,\; f_v + c_e + b_u\bigr) + \kappa ,
\label{eq:elim}
\end{equation}
so $x_e$ is fixed to $0$ whenever $\ell(e) > \overline{z}$; a node is fixed to $0$ once all its
incident edges are eliminated. A complementary demand-reachability test fixes $y_v=0$ when even the
least-demand route to $v$ violates capacity ($2\,g_v - d_v > Q$ for tours, $g_v^{s}+g_v^{t}-d_v > Q$
for paths, with $g$ the minimum cumulative demand to $v$). These run at the root with the warm-start
incumbent and are baked into the formulation.

\paragraph{Domain propagation.}
Inside the tree, propagation runs together with each LP solve (Figure~\ref{fig:flow}) and only once
an incumbent exists. A \emph{sweep} fires on every improvement of $\overline{z}$: it re-applies the
elimination test~\eqref{eq:elim} with the current incumbent, fixing newly inconsistent edges and
nodes. A \emph{chain} propagation fires whenever an edge $(a,i)$ is fixed to one: any other free edge
$e'=(i,j)$ is fixed to $0$ when the cheapest completion forced through both edges,
$f_a + c_{ai} - p_i + c_{ij} + b_j + \kappa$, exceeds $\overline{z}$ (and symmetrically on the
return leg through $a$); with the all-pairs matrix the bound is the minimum over the eight orderings
of the two edges along a depot-rooted tour.

\paragraph{Lagrangian reduced-cost fixing.}
After an optimal LP relaxation of value $z_{\mathrm{LP}}$, the LP reduced costs $\bar{c}$ are used as
labelling weights: the labelling is re-run with edge weights $\bar{c}_e$ and node weights $-\bar{c}_{y_i}$,
giving reduced-cost bounds $\hat{f},\hat{b}$ and the cheapest reduced-cost completion $z_{\mathrm{LR}}$.
A free edge $e$ is fixed to $0$ when
\begin{equation}
z_{\mathrm{LP}} + \Bigl(\min\bigl(\hat{f}_u + \bar{c}_e + \hat{b}_v,\ \hat{f}_v + \bar{c}_e + \hat{b}_u\bigr) - z_{\mathrm{LR}}\Bigr) \;>\; \overline{z},
\label{eq:rcfix}
\end{equation}
the bracket being the extra reduced cost of forcing $e$ in. The pass is guarded to run only when the
Lagrangian gap is sound, $0 \le z_{\mathrm{LP}} - z_{\mathrm{LR}} < (\overline{z}-z_{\mathrm{LP}})$,
since the non-elementary relaxation can otherwise overestimate per-edge penalties. Optionally a node
is fixed to $1$ when forbidding it raises $z_{\mathrm{LR}}$ past $\overline{z}$. The pass supports
four schedules: root-only, on-improvement, periodic, and an \emph{adaptive} schedule that disables
itself after two consecutive improvement rounds yield fewer than five fixings. Algorithm~\ref{alg:rcfix}
summarizes it.

\begin{algorithm}[htbp]
\caption{Reduced-cost fixing at an optimal node LP}
\label{alg:rcfix}
\begin{algorithmic}[1]
\State $z_{\mathrm{LP}} \gets$ LP objective; $\bar{c} \gets$ LP reduced costs; $\overline{z} \gets$ incumbent
\State $\hat f, \hat b \gets$ capacity-aware labelling with weights $\bar{c}$ (forward, backward)
\State $z_{\mathrm{LR}} \gets$ cheapest capacity-feasible completion under $\bar{c}$
\If{$z_{\mathrm{LR}} = \infty$ \textbf{or not} $0 \le z_{\mathrm{LP}}-z_{\mathrm{LR}} < \overline{z}-z_{\mathrm{LP}}$} \Return \Comment{gap unsound}
\EndIf
\For{each free edge $e=(u,v)$}
  \State $\delta \gets \min(\hat f_u + \bar c_e + \hat b_v,\ \hat f_v + \bar c_e + \hat b_u) - z_{\mathrm{LR}}$
  \If{$z_{\mathrm{LP}} + \delta > \overline{z}$} fix $x_e \gets 0$ \EndIf
\EndFor
\State fix any node whose incident edges are all fixed; \textbf{optionally} fix $y_i \gets 1$ per~\eqref{eq:rcfix}
\end{algorithmic}
\end{algorithm}

These components are individually switchable, which is what makes the ablation in
Section~\ref{sec:study} possible.

\subsection{Primal heuristics and branching}
\label{sec:heur}

The cut and bound machinery is supported by primal heuristics and a problem-specific branching
scheme. We describe them briefly; they are enabled by default in every experiment but are not the
subject of the ablation.

\paragraph{Primal heuristics.}
The warm-start heuristic (Algorithm~\ref{alg:heur}) first enumerates all one- and two-customer seed
routes through the terminals, ranks them by objective, and extends each by cheapest-insertion: a
customer $c$ is inserted into the gap $(a,b)$ minimizing $c_{ac} + c_{cb} - c_{ab}$ subject to
capacity. Each construction is improved by iterated local search whose neighbourhoods are 2-opt,
Or-opt (relocation of chains of length $1$--$3$), a profit-aware swap of a visited for an unvisited
customer, and profit-aware drop/add moves: a customer is dropped when its incident edge saving
exceeds its profit, and added when the reverse holds. Rounds of a deterministic perturbation kick
(swap or relocate) alternate with local search, accepting only improving solutions. The whole
heuristic is deterministic (kicks come from a hash of the iteration indices, not a clock), so runs are
reproducible.

During the search, an LP-guided heuristic runs in the branch-and-bound callback (rate-limited to once
per $200$ nodes) on the graph reduced to the fractional support. Three reductions seed it from
$(\bar x,\bar y)$: a \emph{threshold} mode keeps nodes with $\bar y_i > 0.5$ and edges with
$\bar x_e > 0.1$; a \emph{RINS} mode keeps the edges and nodes on which the incumbent and the LP agree
or which the LP leaves fractional; and a \emph{neighbourhood} mode seeds from edges with $\bar x_e >
0.3$ and adds all edges among the seeded nodes. Construction (in LP-weighted order, customers sorted
by $\bar y_i$) and the same iterated local search then run on the reduced subgraph; the best tour, if
it improves the incumbent and passes a subtour-feasibility check, is handed to the solver. Good
incumbents matter beyond their own objective value: every bound-based elimination and propagation in
Section~\ref{sec:fixing} is driven by the gap to the incumbent.

\paragraph{Branching.}
By default the solver uses HiGHS's pseudocost variable branching. The implementation also offers
\emph{hyperplane branching}, which branches on a linear expression $a^{\!\top} y$ in the node-selection
variables. At a fractional point the two children are the rounded half-spaces
\begin{equation}
a^{\!\top} y \;\le\; \bigl\lfloor a^{\!\top}\bar y \bigr\rfloor
\qquad\text{and}\qquad
a^{\!\top} y \;\ge\; \bigl\lceil a^{\!\top}\bar y \bigr\rceil ,
\label{eq:hyper}
\end{equation}
instantiated in four families: a \emph{pair} $a^{\!\top} y = y_i + y_j$ for two near-cost customers
(a Ryan--Foster-style same/different split, here on the node sum); a demand-seeded \emph{cluster}
$\sum_{i \in C} d_i y_i$; the \emph{total selected demand} $\sum_i d_i y_i$; and the \emph{cardinality}
$\sum_i y_i$. Candidates are scored in three tiers: a dual-based prescore
$\bigl(\sum_j |a_j \bar{c}_{y_j}|\bigr)\cdot\min(\phi,1-\phi)$ (with $\phi$ the fractionality of
$a^{\!\top}\bar y$) keeps the top three; reliable online pseudocosts (both directions sampled at least
four times) score those directly; and strong branching resolves the rest (this third tier is off by
default). Hyperplane branching overrides the variable branch only when its best score beats the
variable score. The experiments use the default (variable) branching configuration throughout, so the
hyperplane machinery is described for completeness, not exercised in Section~\ref{sec:study}.

\begin{algorithm}[htbp]
\caption{Warm-start construction and iterated local search}
\label{alg:heur}
\begin{algorithmic}[1]
\State $\mathcal{P} \gets$ all $1$- and $2$-customer seed routes through the terminals, capacity-feasible, ranked by objective
\State $z^\star \gets \infty$
\For{each seed $\sigma$ in $\mathcal{P}$ (in parallel)}
  \State extend $\sigma$ by cheapest-insertion: repeatedly add the customer/gap minimizing $c_{ac}+c_{cb}-c_{ab}$ within capacity
  \State $\sigma \gets \textsc{LocalSearch}(\sigma)$ \Comment{2-opt, Or-opt(1--3), swap, drop/add}
  \For{a fixed number of rounds}
    \State $\sigma' \gets \textsc{LocalSearch}(\textsc{Kick}(\sigma))$; \textbf{if} $\mathrm{obj}(\sigma') < \mathrm{obj}(\sigma)$ \textbf{then} $\sigma \gets \sigma'$
  \EndFor
  \State $z^\star \gets \min(z^\star, \mathrm{obj}(\sigma))$
\EndFor
\State \Return the best tour found
\end{algorithmic}
\end{algorithm}

\section{Computational study}
\label{sec:study}

\paragraph{Setup.}
All experiments are run on an AMD Ryzen~9 3950X (16 cores, 32 threads, 124\,GB RAM) under Ubuntu
24.04, with the solver compiled in Release mode with GCC~14.2 (C++23) against HiGHS
v1.15.0~\cite{huangfu2018highs}. The solver and every benchmark script are archived as open
software~\cite{cptp_software}. The head-to-head comparison with dynamic programming
(Section~\ref{sec:vsdp}) uses a per-instance wall-clock limit of $3600$ seconds; the ablation
(Sections~\ref{sec:ablation}--\ref{sec:spi-result}), which multiplies runs across configurations,
uses a reduced limit of $300$ seconds, stated again where used. Unless noted otherwise, separation,
heuristics, and branching use their default settings; each ablation configuration is reported
together with the non-default options it sets, and every result table is regenerated directly from
the committed run logs.

\paragraph{What we do not report.}
The solver models a single node-demand capacity, so multi-resource instances, time-windowed
variants, and arc-resource formulations are out of scope (Section~\ref{sec:instances}). We report
wall-clock runtime, branch-and-bound nodes, and per-family cut counts; we do not study parallel
scaling, and all runs use the solver's default thread count on the machine above.

\subsection{Instances}
\label{sec:instances}
We use two instance sets, both column-generation pricing subproblems that arise in vehicle-routing,
taken from the \texttt{cptp} repository and used in the capacitated-profitable-tour study of Jepsen
et al.~\cite{jepsen2014cptp}. The SPPRCLIB set (45 instances) comprises the elementary
shortest-path-with-capacity instances of Jepsen, Petersen and Spoorendonk~\cite{jepsen2008spprc};
the Roberti Set~3 set (31 instances) was generated for that study by R.~Roberti using the pricing
code of Baldacci, Mingozzi and Roberti~\cite{baldacci2011ng}. Both carry the negative reduced costs
of routing pricing, so elementarity binds and they match the single-capacity model of
Section~\ref{sec:model} directly. Table~\ref{tbl:instancefamily} summarizes their sizes.

\begin{table}[htbp]
\centering
\caption{Benchmark instance families: number of instances and the range of node count, edge count, and capacity within each family.}
\label{tbl:instancefamily}
\vspace{0.5em}
\begin{tabular}{lrrrr}
\toprule
Family & \# & Nodes & Edges & Capacity \\
\midrule
SPPRCLIB & 45 & 45--262 & 990--34191 & 70--500 \\
Roberti Set 3 & 31 & 45--200 & 990--19900 & 100--30000 \\
\bottomrule
\end{tabular}
\end{table}

We considered the classical resource-constrained shortest-path instances of Beasley and
Christofides~\cite{beasley1989rcsp} (OR-Library) and the negative-cost classes of
García~\cite{garcia2009rcsp} as external comparators, but they do not fit the model of
Section~\ref{sec:model}. Their resource is consumed on \emph{arcs} over a \emph{directed} graph,
whereas the model here is an undirected formulation with a single \emph{node}-demand capacity; arc
consumption cannot be re-expressed as node consumption, the non-reducibility already noted by Jepsen
et al.~\cite{jepsen2014cptp}. García's random classes that most resemble the profitable-tour regime
are moreover multi-resource ($R \in \{5,10\}$) and not publicly available. Importing these
benchmarks faithfully would require a different model (per-arc resources, directed arcs), which we
do not pursue here. The solver is open, however, and this is a localized change (to the model and
the cut separators) that we would encourage others to make should such tests be of interest.

\subsection{Branch-and-cut versus dynamic programming}
\label{sec:vsdp}
We compare the branch-and-cut solver, in its default configuration (GSEC, RCI, and multistar cuts;
preprocessing, propagation, and reduced-cost fixing enabled), against a dynamic-programming /
labelling reference (PathWyse~\cite{salani2024pathwyse}) on the $76$ instances of the SPPRCLIB and
Roberti sets. We use PathWyse rather than the author's companion bucket-graph SPPRC
library~\cite{bgspprc2026} as the labelling reference: obtaining exact elementarity from a
bucket-graph labelling requires incrementally enlarging the node neighbourhoods, a scheme PathWyse
implements but the bucket-graph library, which fixes the neighbourhood size a priori, does not, so
it is not directly applicable to the elementary problem solved here.
Both solvers prove optimality on almost the entire benchmark: PathWyse on $71$ of $76$
instances and the branch-and-cut solver on $69$ of $76$. Every instance either is solved to proven
optimality or hits the time limit; the branch-and-cut solver closes the optimality gap to zero on
all $69$ instances it solves (its objectives match PathWyse on every instance both solve), so the
differences below are differences in runtime, not in solution quality.

\begin{figure}[htbp]
\centering
\includegraphics[width=0.7\textwidth]{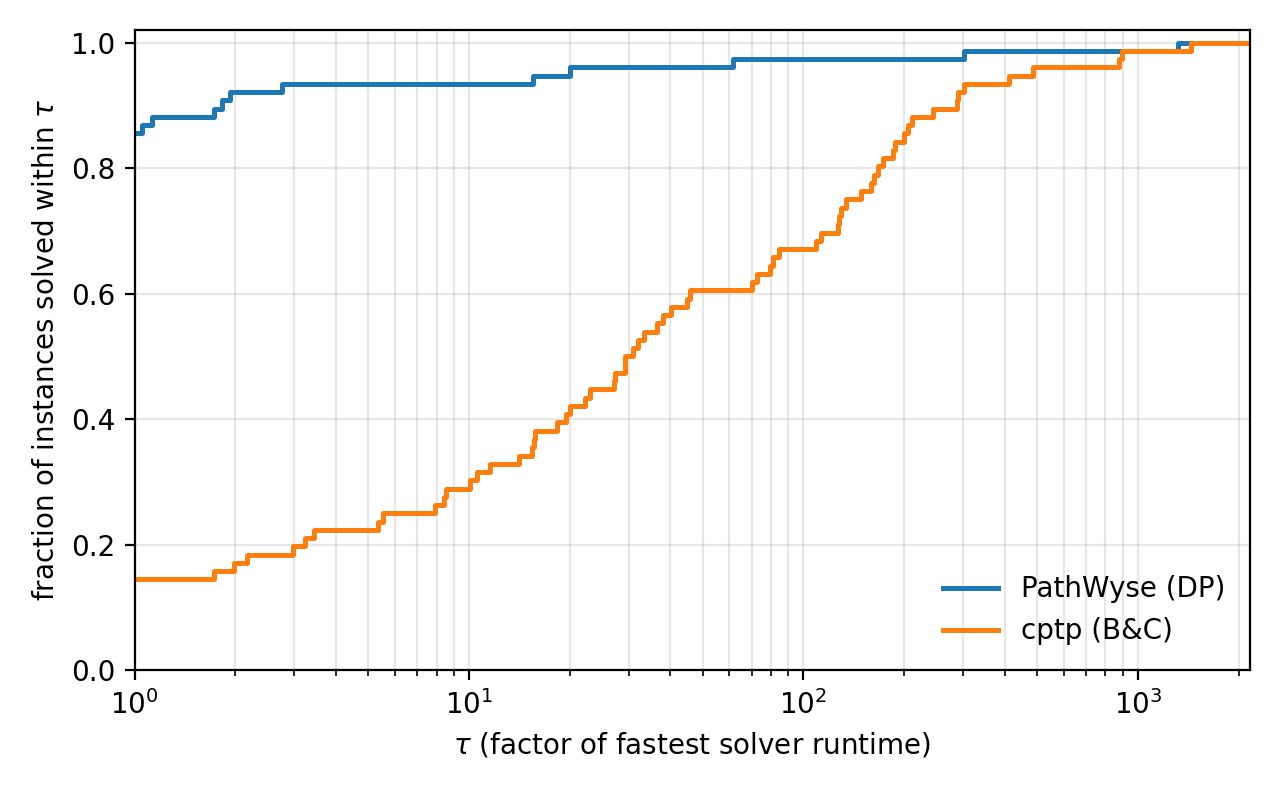}
\caption{Performance profile (Dolan--Moré) of runtime over the $76$ common instances:
$\rho(\tau)$ is the fraction of instances solved within a factor $\tau$ of the fastest solver.}
\label{fig:perfprofile}
\end{figure}

\begin{figure}[htbp]
\centering
\includegraphics[width=0.7\textwidth]{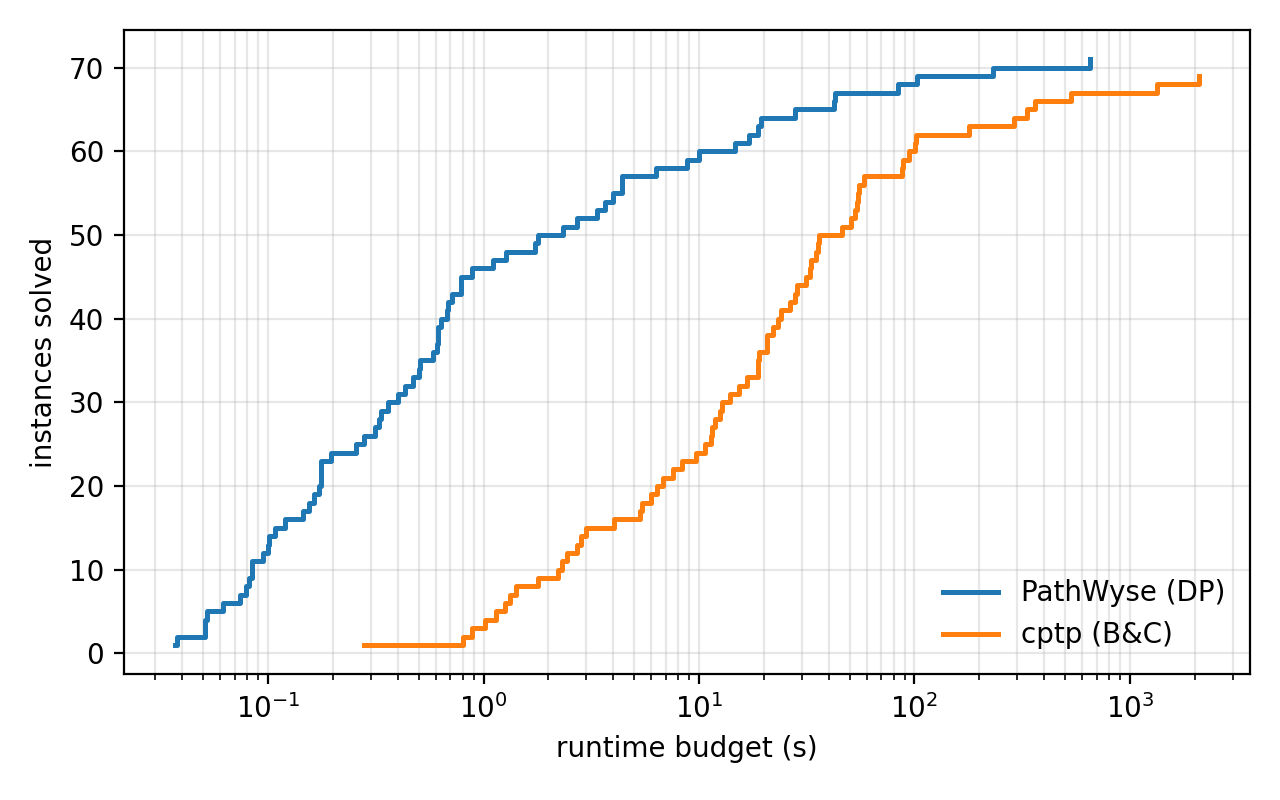}
\caption{Cactus plot: number of instances solved within a given wall-clock budget (log scale).}
\label{fig:cactusruntime}
\end{figure}

On the instances both solve, labelling is the faster method by a clear margin. The performance
profile (Figure~\ref{fig:perfprofile}) and the cactus plot (Figure~\ref{fig:cactusruntime}) show
PathWyse ahead across the runtime range, and the shifted geometric mean of runtime
(Table~\ref{tbl:sgmstats}) is $3.7$\,s for PathWyse against $31.1$\,s for branch-and-cut (a factor
of $8.5$), with a larger gap at the median ($0.6$\,s against $21.4$\,s). On the easy,
elementarity-only instances, where the resource does not bind and few labels survive dominance,
there is no reason to prefer a cutting-plane method, and we do not.

The value of branch-and-cut shows instead at the hard tail, where the two methods are
\emph{complementary}. The five instances on which PathWyse exhausts the one-hour limit are all
solved by branch-and-cut, several of them quickly: \texttt{B-n57-k7-20} in $2.7$\,s,
\texttt{M-n200-k16-143} in $11.9$\,s, \texttt{G-n262-k25-316} in $58$\,s, \texttt{M-n121-k7-260} in
$180$\,s, and \texttt{M-n121-k7\_a} in about $35$ minutes, still within the hour. Conversely, the
seven instances on which branch-and-cut times out, four of them from the \texttt{M-n200-k16/k17}
family, are solved by PathWyse in seconds to minutes. Neither method dominates; each closes the
other's open instances, echoing the
complementarity reported by Jepsen et al.~\cite{jepsen2014cptp}.

\begin{table}[htbp]
\centering
\caption{Shifted geometric mean (SGM, shift $s=1$\,s) and median of runtime over the 76 common instances; unsolved runs enter at the time limit. ``Scaled'' normalizes the SGM to the fastest solver.}
\label{tbl:sgmstats}
\vspace{0.5em}
\begin{tabular}{lrrrr}
\toprule
Solver & Solved & SGM time (s) & Median time (s) & Scaled \\
\midrule
PathWyse (DP) & 71/76 & 3.66 & 0.61 & 1.00 \\
cptp (B\&C) & 69/76 & 31.10 & 21.42 & 8.49 \\
\bottomrule
\end{tabular}
\end{table}

\subsection{Gating ablation: which cuts and fixing pay off}
\label{sec:ablation}
We isolate the marginal value of each component by enabling them incrementally, starting from
connectivity (GSEC) cuts only and adding, in turn, the capacity-class cuts (RCI and multistar),
reduced-cost fixing and bound-based propagation, the comb and rounded-GLM cuts, the all-pairs
bounds, and finally the SPI separator. Every configuration is run on all $76$ instances under a
$300$-second limit; Table~\ref{tbl:ablation} reports the result.

The capacity-class cuts account for essentially the entire improvement. With GSEC alone the solver
proves optimality on $52$ of $76$ instances; adding RCI and multistar raises this to $64$, cuts the
shifted geometric mean of the branch-and-bound tree by more than an order of magnitude (from $2290$
to $214$ nodes) and the runtime by a factor of $2.4$ (from $54.7$ to $22.8$\,s), and, because the
tree shrinks so much, reduces the total number of cuts added by nearly tenfold. This is the
central positive result: on capacity-binding instances the capacity-class inequalities are
decisive, exactly the regime in which a cutting-plane approach earns its place against labelling.

The remaining components add nothing measurable on this benchmark. Reduced-cost fixing and
bound-based propagation leave the solved count, node count, and runtime essentially unchanged, and
this is robust to where they enter the ladder (Table~\ref{tbl:fixrobust}): switching them on nudges
the shifted geometric mean from $22.8$ to $23.1$\,s on top of the capacity cuts and from $22.9$ to
$23.2$\,s on top of the comb and rounded-GLM cuts, a negligible slowdown, with the search tree
unchanged at both positions. The comb and
rounded-GLM cuts likewise do not help. We keep
them available (they are inexpensive and may matter on instances outside this set) but report
plainly that they are not what makes the method work here.

\begin{table}[htbp]
\centering
\caption{Gating ablation over the 76 SPPRCLIB and Roberti instances (per-instance limit 300\,s). Each row adds one component to the row above. SGM is the shifted geometric mean (shift $s=1$\,s for time); unsolved runs enter at the time limit. Capacity-class cuts account for essentially the entire improvement; reduced-cost fixing, comb/RGLM, and SPI add nothing, and the all-pairs bounds add a slight overhead.}
\label{tbl:ablation}
\vspace{0.5em}
\resizebox{\ifdim\width>\linewidth\linewidth\else\width\fi}{!}{%
\begin{tabular}{lrrrr}
\toprule
Configuration & Solved & SGM time (s) & SGM nodes & Mean cuts \\
\midrule
GSEC only & 52/76 & 54.7 & 2290 & 25875 \\
\quad + capacity cuts (RCI, multistar) & 64/76 & 22.8 & 214 & 2692 \\
\quad + reduced-cost fixing & 64/76 & 23.1 & 217 & 2753 \\
\quad + comb, rounded GLM & 64/76 & 23.2 & 218 & 2910 \\
\quad + all-pairs bounds & 63/76 & 24.6 & 215 & 2714 \\
\quad + SPI & 63/76 & 24.5 & 215 & 2713 \\
\bottomrule
\end{tabular}
}
\end{table}

\begin{table}[htbp]
\centering
\caption{Reduced-cost fixing and bound-based propagation are inert regardless of where they enter the ladder: shifted geometric mean runtime and search-tree size over the 76 instances (300\,s limit) with the fixing/propagation bundle off vs.\ on, at two positions in the ladder. Turning it on is a negligible slowdown and shrinks the tree at neither position.}
\label{tbl:fixrobust}
\vspace{0.5em}
\begin{tabular}{llrrr}
\toprule
Cut set & Fixing/prop.\ & Solved & SGM time (s) & SGM nodes \\
\midrule
GSEC + capacity cuts & off & 64/76 & 22.8 & 214 \\
GSEC + capacity cuts & on & 64/76 & 23.1 & 217 \\
\quad + comb, rounded GLM & off & 64/76 & 22.9 & 214 \\
\quad + comb, rounded GLM & on & 64/76 & 23.2 & 218 \\
\bottomrule
\end{tabular}
\end{table}

Table~\ref{tbl:cutactivity} shows the mechanism behind these numbers. Without capacity-class cuts
the separator floods the cut pool, almost $26{,}000$ GSEC cuts per instance on average, as
connectivity inequalities alone struggle to tighten the relaxation; adding the rounded-capacity and
multistar cuts replaces roughly an order of magnitude of that volume with about $1{,}300$
capacity-class cuts. The components that do not help are precisely the ones that almost never
separate anything: comb inequalities find \emph{no} violated cut on any instance, RGLM separates
barely two cuts per instance, and SPI (Section~\ref{sec:spi-result}) finds none. Their failure here
is one of separation, not of strength. Table~\ref{tbl:cuttime} gives the matching per-separator
times: the off-by-default families cost almost nothing to run (comb and SPI under $0.01$\,s per
instance, rounded GLM about $0.3$\,s), so leaving them enabled would be nearly free.

\begin{table}[htbp]
\centering
\caption{Mean number of cuts added per instance, by separator and configuration (76 instances, $300$\,s limit). Capacity-class cuts replace an order of magnitude of GSEC cuts; comb, RGLM, and SPI separate almost nothing, which is why they do not affect performance.}
\label{tbl:cutactivity}
\vspace{0.5em}
\begin{tabular}{lrrrrrr}
\toprule
Configuration & GSEC & RCI & Multistar & Comb & RGLM & SPI \\
\midrule
GSEC only & 25681 & -- & -- & -- & -- & -- \\
+ capacity cuts & 1336 & 8.4 & 1265 & -- & -- & -- \\
+ comb, RGLM & 1475 & 8.5 & 1340 & -- & 1.6 & -- \\
+ all-pairs, SPI & 1348 & 8.4 & 1275 & -- & -- & -- \\
\bottomrule
\end{tabular}
\end{table}

\begin{table}[htbp]
\centering
\caption{Mean separation time (s) per instance, by separator and configuration (76 instances, $300$\,s limit; ``--'' where a separator is disabled). Separators run in parallel, so these are per-family costs, not a sum. The off-by-default families are cheap: comb and SPI cost under $0.01$\,s and rounded GLM about $0.3$\,s, negligible against the LP and branch-and-bound that dominate the runtime (Table~\ref{tbl:componentcost}).}
\label{tbl:cuttime}
\vspace{0.5em}
\begin{tabular}{lrrrrrr}
\toprule
Configuration & GSEC & RCI & Multistar & Comb & RGLM & SPI \\
\midrule
GSEC only & 1.07 & -- & -- & -- & -- & -- \\
+ capacity cuts & 0.33 & 0.36 & 0.32 & -- & -- & -- \\
+ comb, RGLM & 0.34 & 0.39 & 0.34 & $<$0.01 & 0.27 & -- \\
+ all-pairs, SPI & 0.33 & 0.37 & 0.32 & -- & -- & $<$0.01 \\
\bottomrule
\end{tabular}
\end{table}

Finally, the components are cheap in wall-clock terms, which is why the default keeps the useful
cuts on and the rest available rather than removed. Table~\ref{tbl:componentcost} breaks the default
configuration's runtime down by component: the LP relaxation and branch-and-bound account for about
$98\%$ of it, the primal heuristic for a further $1.4\%$, and all cut separation, domain
propagation, and reduced-cost fixing together for about one percent. The cut families that are off
by default are similarly inexpensive to separate (Table~\ref{tbl:cuttime}); they are off only
because they yield no measurable benefit here (Table~\ref{tbl:ablation}), not because they are
expensive, so we keep them available for instances where they might pay off. Reduced-cost fixing tells the same story from the other side: it fixes on
the order of $65$ variables per instance yet does not shorten the search, because the capacity cuts
already keep the tree small, and it is retained because it is nearly free. The default is therefore
the productive GSEC, RCI, and multistar cuts together with the near-free propagation and
reduced-cost fixing, and nothing else.

\begin{table}[htbp]
\centering
\caption{Where the default configuration (GSEC, RCI, and multistar cuts with domain propagation and reduced-cost fixing) spends its runtime: mean wall-clock per instance over the 76 instances (3600\,s limit), by component. The LP relaxation and branch-and-bound dominate; all separation, propagation, and fixing together are about one percent. Separators run in parallel, so the separation figure is the wall-clock phase span, not the sum of per-family times.}
\label{tbl:componentcost}
\vspace{0.5em}
\begin{tabular}{llrr}
\toprule
Component & Activity & Time (s) & Share \\
\midrule
LP relaxation and branch-and-bound &  & 408.43 & 97.6\% \\
Primal heuristic &  & 5.96 & 1.4\% \\
Cut separation (GSEC, RCI, multistar) & $\approx$~12467 cuts & 3.86 & 0.9\% \\
Domain propagation and reduced-cost fixing & $\approx$~66 fixings & 0.34 & 0.1\% \\
\midrule
Total & & 418.6 & 100\% \\
\bottomrule
\end{tabular}
\end{table}

\subsection{A negative result on SPI}
\label{sec:spi-result}
The shortest-path-incompatibility separator is the one component for which we expected a payoff on
the negative-cost, capacity-binding instances, and it delivers none. The all-pairs bounds that SPI
depends on already cost a little on their own: enabling them lifts the shifted geometric mean from
$23.2$ to $24.6$\,s and loses one instance to the time limit (the solved count falls from $64$ to
$63$ of $76$). Adding the SPI cut on top of them changes nothing measurable: $63$ of $76$ solved,
SGM $24.5$\,s, indistinguishable from the $24.6$\,s without it. The reason is direct: \emph{the SPI
separator finds zero violated cuts on all $76$ instances}. The shortest-path/Held--Karp incompatibility bound is never tight enough, relative to
the incumbent, to exclude a pair or set of nodes that the LP relaxation has not already excluded.
We report this as a negative result so that the cut (a variant of García's node-precedence
inequalities~\cite{garcia2009rcsp}, here tightened with a Held--Karp bound) is not mistaken for a
useful addition on this problem class.

\section{Conclusion}
\label{sec:conclusion}
We have presented an open, reproducible branch-and-cut for the capacitated profitable tour problem
and its open $s$--$t$ path variant, built end-to-end on a fully open MIP stack, and used it for a
component study. Three findings stand out. First, on a common modern stack the
branch-and-cut reproduces the long-standing picture of being complementary to dynamic programming:
labelling is faster on the easy instances, but the two methods close each other's hard instances, so
neither dominates. Second, that complementarity is driven almost entirely by the capacity-class
cuts (GSEC alone solves $52$ of $76$ instances, adding RCI and multistar raises this to $64$ and
shrinks the search tree more than tenfold), while comb and rounded-GLM cuts, reduced-cost fixing,
and bound-based propagation add nothing measurable here. Third, the shortest-path-incompatibility
cut, a variant of García's node-precedence inequalities tightened with a Held--Karp bound, is a dead
cut on this problem class: it finds no violated inequality on any of the $76$ instances.

We make no claim of a new method: the branch-and-cut is that of Jepsen et
al.~\cite{jepsen2014cptp}, whose CPTP formulation and cut families it re-implements; the one
component drawn from García~\cite{garcia2009rcsp}, the SPI cut, contributes nothing here. The
contribution is the open, reproducible artifact~\cite{cptp_software}, a rerunnable baseline others
can build on, together with the decomposition of which components pay off and where its
time goes, including the
negative results, which we hope save others from re-deriving them. Natural next steps are to extend
the formulation to per-arc resources and directed graphs (which would bring the classical
resource-constrained shortest-path benchmarks into range) and to test whether the capacity-class
cuts retain their decisive role on larger and more strongly capacity-binding instances.

\section*{Code and data availability}
The solver is open source (MIT) at \url{https://github.com/spoorendonk/cptp} and archived on
Zenodo~\cite{cptp_software}. Its \texttt{benchmarks/} directory contains the run, ablation, and
PathWyse-reference scripts, the reference optima, and the committed result files
(\texttt{cptp.csv}, \texttt{pathwyse.csv}, \texttt{ablation.csv}) from which every table and figure
in this paper is regenerated. The SPPRCLIB and Roberti pricing instances are distributed with that
repository; PathWyse~\cite{salani2024pathwyse} is cloned from its upstream repository by the setup
script.

The solver is driven from the command line or from Python. Building it and reproducing a result,
the configuration toggles used in the ablation, and the full benchmark and ablation runs are a
handful of commands:

\begin{lstlisting}
# Build (HiGHS is fetched and patched automatically)
cmake -B build -DCMAKE_BUILD_TYPE=Release && cmake --build build -j

# Solve one instance (tour); --source/--target gives the open s-t path variant
./build/cptp-solve benchmarks/instances/spprclib/B-n45-k6-54.sppcc
./build/cptp-solve instance.txt --source 0 --target 12

# Toggle cut families, bounds, and fixing (the ablation configurations)
./build/cptp-solve <instance> --enable_comb true --enable_rglm true \
    --all_pairs_bounds true --enable_spi true --rc_fixing adaptive

# Regenerate the result CSVs behind every table and figure
./benchmarks/run_benchmarks.sh --time-limit 3600   # cptp.csv
./benchmarks/run_ablation.sh   --time-limit 300    # ablation.csv
\end{lstlisting}

\noindent The same is available through the Python bindings:

\begin{lstlisting}[language=Python]
import cptp
prob   = cptp.load("instance.txt")            # .txt, .vrp, or .sppcc
model  = cptp.Model(); model.set_problem(prob)
result = model.solve([("all_pairs_bounds", "true"),
                      ("enable_spi", "true"),
                      ("time_limit", "60")])
print(result.objective, result.tour)
\end{lstlisting}

\section*{Acknowledgements}
We thank R.~Roberti for the Set~3 pricing instances, and the developers of
HiGHS~\cite{huangfu2018highs} and PathWyse~\cite{salani2024pathwyse} for the open-source solvers on
which this study is built.

\appendix
\section{Per-instance results}
\label{sec:appendix-results}
Table~\ref{tbl:instancedetail} lists the characteristics of every instance. Table~\ref{tbl:runtime}
gives the per-instance objective and runtime for branch-and-cut and dynamic programming
(Section~\ref{sec:vsdp}), and Table~\ref{tbl:ablationinstances} the per-instance runtime under each
ablation configuration (Section~\ref{sec:ablation}).

\begingroup\footnotesize
\setlength{\tabcolsep}{4pt}
\begin{longtable}{lrrr}
\caption{Per-instance characteristics: node count, edge count, and capacity.}\label{tbl:instancedetail}\\
\toprule
Instance & Nodes & Edges & Capacity \\
\midrule
\endfirsthead
\multicolumn{4}{c}{\tablename\ \thetable\ -- continued from previous page}\\
\toprule
Instance & Nodes & Edges & Capacity \\
\midrule
\endhead
\midrule\multicolumn{4}{r}{\textit{continued on next page}}\\
\endfoot
\bottomrule
\endlastfoot
\midrule
\multicolumn{4}{l}{\textbf{SPPRCLIB}} \\
\texttt{B-n45-k6-54} & 45 & 990 & 100 \\
\texttt{B-n50-k8-40} & 50 & 1225 & 100 \\
\texttt{P-n50-k10-24} & 50 & 1225 & 100 \\
\texttt{P-n50-k7-92} & 50 & 1225 & 150 \\
\texttt{P-n50-k8-19} & 50 & 1225 & 120 \\
\texttt{P-n51-k10-30} & 51 & 1275 & 80 \\
\texttt{B-n52-k7-15} & 52 & 1326 & 100 \\
\texttt{A-n54-k7-149} & 54 & 1431 & 100 \\
\texttt{P-n55-k10-44} & 55 & 1485 & 115 \\
\texttt{P-n55-k15-88} & 55 & 1485 & 70 \\
\texttt{P-n55-k7-116} & 55 & 1485 & 170 \\
\texttt{P-n55-k8-260} & 55 & 1485 & 160 \\
\texttt{B-n57-k7-20} & 57 & 1596 & 100 \\
\texttt{A-n60-k9-57} & 60 & 1770 & 100 \\
\texttt{P-n60-k10-24} & 60 & 1770 & 120 \\
\texttt{P-n60-k15-8} & 60 & 1770 & 80 \\
\texttt{A-n61-k9-80} & 61 & 1830 & 100 \\
\texttt{A-n62-k8-99} & 62 & 1891 & 100 \\
\texttt{A-n63-k10-44} & 63 & 1953 & 100 \\
\texttt{A-n63-k9-157} & 63 & 1953 & 100 \\
\texttt{A-n64-k9-45} & 64 & 2016 & 100 \\
\texttt{A-n65-k9-10} & 65 & 2080 & 100 \\
\texttt{P-n65-k10-102} & 65 & 2080 & 130 \\
\texttt{B-n66-k9-50} & 66 & 2145 & 100 \\
\texttt{B-n67-k10-26} & 67 & 2211 & 100 \\
\texttt{B-n68-k9-65} & 68 & 2278 & 100 \\
\texttt{A-n69-k9-42} & 69 & 2346 & 100 \\
\texttt{P-n70-k10-12} & 70 & 2415 & 135 \\
\texttt{E-n76-k10-72} & 76 & 2850 & 140 \\
\texttt{E-n76-k14-102} & 76 & 2850 & 100 \\
\texttt{E-n76-k15-40} & 76 & 2850 & 100 \\
\texttt{E-n76-k7-44} & 76 & 2850 & 220 \\
\texttt{P-n76-k4-41} & 76 & 2850 & 350 \\
\texttt{P-n76-k5-16} & 76 & 2850 & 280 \\
\texttt{B-n78-k10-70} & 78 & 3003 & 100 \\
\texttt{A-n80-k10-14} & 80 & 3160 & 100 \\
\texttt{E-n101-k14-158} & 101 & 5050 & 112 \\
\texttt{E-n101-k8-291} & 101 & 5050 & 200 \\
\texttt{M-n101-k10-97} & 101 & 5050 & 200 \\
\texttt{P-n101-k4-174} & 101 & 5050 & 400 \\
\texttt{M-n121-k7-260} & 121 & 7260 & 200 \\
\texttt{M-n151-k12-15} & 151 & 11325 & 200 \\
\texttt{M-n200-k16-143} & 200 & 19900 & 200 \\
\texttt{M-n200-k17-12} & 200 & 19900 & 200 \\
\texttt{G-n262-k25-316} & 262 & 34191 & 500 \\
\midrule
\multicolumn{4}{l}{\textbf{Roberti Set 3}} \\
\texttt{F-n45-k4\_a} & 45 & 990 & 2010 \\
\texttt{P-n70-k10\_a} & 70 & 2415 & 135 \\
\texttt{P-n70-k10\_b} & 70 & 2415 & 135 \\
\texttt{F-n72-k4\_a} & 72 & 2556 & 30000 \\
\texttt{E-n76-k10\_a} & 76 & 2850 & 140 \\
\texttt{E-n76-k10\_b} & 76 & 2850 & 140 \\
\texttt{E-n76-k14\_a} & 76 & 2850 & 100 \\
\texttt{E-n76-k14\_b} & 76 & 2850 & 100 \\
\texttt{E-n76-k7\_a} & 76 & 2850 & 220 \\
\texttt{E-n76-k7\_b} & 76 & 2850 & 220 \\
\texttt{E-n76-k8\_a} & 76 & 2850 & 180 \\
\texttt{E-n76-k8\_b} & 76 & 2850 & 180 \\
\texttt{P-n76-k4\_a} & 76 & 2850 & 350 \\
\texttt{P-n76-k4\_b} & 76 & 2850 & 350 \\
\texttt{P-n76-k5\_a} & 76 & 2850 & 280 \\
\texttt{P-n76-k5\_b} & 76 & 2850 & 280 \\
\texttt{E-n101-k14\_a} & 101 & 5050 & 112 \\
\texttt{E-n101-k14\_b} & 101 & 5050 & 112 \\
\texttt{E-n101-k8\_a} & 101 & 5050 & 200 \\
\texttt{E-n101-k8\_b} & 101 & 5050 & 200 \\
\texttt{P-n101-k4\_a} & 101 & 5050 & 400 \\
\texttt{P-n101-k4\_b} & 101 & 5050 & 400 \\
\texttt{M-n121-k7\_a} & 121 & 7260 & 200 \\
\texttt{M-n121-k7\_b} & 121 & 7260 & 200 \\
\texttt{F-n135-k7\_a} & 135 & 9045 & 2210 \\
\texttt{M-n151-k12\_a} & 151 & 11325 & 200 \\
\texttt{M-n151-k12\_b} & 151 & 11325 & 200 \\
\texttt{M-n200-k16\_a} & 200 & 19900 & 200 \\
\texttt{M-n200-k16\_b} & 200 & 19900 & 200 \\
\texttt{M-n200-k17\_a} & 200 & 19900 & 200 \\
\texttt{M-n200-k17\_b} & 200 & 19900 & 200 \\
\end{longtable}
\endgroup

\begingroup\footnotesize
\setlength{\tabcolsep}{4pt}
\begin{longtable}{lrrr}
\caption{Per-instance results: the optimal objective (from whichever solver proves optimality) and runtime (s) for branch-and-cut (cptp, default configuration) and dynamic programming (PathWyse), at a $3600$\,s limit. $\dagger$ marks a timeout.}\label{tbl:runtime}\\
\toprule
Instance & Objective & cptp (s) & PathWyse (s) \\
\midrule
\endfirsthead
\multicolumn{4}{c}{\tablename\ \thetable\ -- continued from previous page}\\
\toprule
Instance & Objective & cptp (s) & PathWyse (s) \\
\midrule
\endhead
\midrule\multicolumn{4}{r}{\textit{continued on next page}}\\
\endfoot
\bottomrule
\endlastfoot
\midrule
\multicolumn{4}{l}{\textbf{SPPRCLIB}} \\
\texttt{P-n70-k10-12} & -70317 & 1.15 & 0.14 \\
\texttt{P-n76-k5-16} & -107633 & 1.33 & 3.68 \\
\texttt{P-n50-k8-19} & -83307 & 1.41 & 0.71 \\
\texttt{P-n60-k10-24} & -15001 & 2.22 & 0.10 \\
\texttt{P-n55-k7-116} & -17824 & 2.31 & 0.08 \\
\texttt{B-n57-k7-20} & -867154 & 2.73 & 3600$^\dagger$ \\
\texttt{B-n67-k10-26} & -21924 & 2.83 & 0.28 \\
\texttt{B-n52-k7-15} & -74998 & 2.99 & 1.73 \\
\texttt{A-n69-k9-42} & -43290 & 4.04 & 0.18 \\
\texttt{A-n63-k10-44} & -32561 & 5.33 & 0.20 \\
\texttt{B-n45-k6-54} & -74278 & 5.44 & 9.97 \\
\texttt{A-n65-k9-10} & -42835 & 6.36 & 0.33 \\
\texttt{E-n76-k7-44} & -22214 & 6.84 & 0.43 \\
\texttt{P-n76-k4-41} & -88276 & 8.36 & 8.80 \\
\texttt{A-n63-k9-157} & -24189 & 9.73 & 0.12 \\
\texttt{P-n50-k10-24} & -2965 & 10.69 & 0.04 \\
\texttt{P-n55-k8-260} & -3573 & 11.36 & 0.10 \\
\texttt{A-n61-k9-80} & -23549 & 11.48 & 0.31 \\
\texttt{M-n200-k16-143} & -198792 & 11.88 & 3600$^\dagger$ \\
\texttt{A-n54-k7-149} & -12492 & 12.53 & 0.68 \\
\texttt{E-n76-k10-72} & -25241 & 12.85 & 0.40 \\
\texttt{P-n50-k7-92} & -2 & 13.98 & 0.18 \\
\texttt{B-n50-k8-40} & -12832 & 15.31 & 0.34 \\
\texttt{P-n101-k4-174} & -17702 & 16.64 & 18.80 \\
\texttt{M-n101-k10-97} & -32628 & 18.78 & 6.32 \\
\texttt{P-n51-k10-30} & -2 & 20.64 & 0.16 \\
\texttt{P-n55-k10-44} & -1090 & 20.78 & 0.15 \\
\texttt{B-n78-k10-70} & -44333 & 23.17 & 2.71 \\
\texttt{P-n60-k15-8} & -534 & 24.11 & 0.08 \\
\texttt{A-n64-k9-45} & -50550 & 26.43 & 0.79 \\
\texttt{E-n101-k8-291} & -4266 & 28.49 & 0.63 \\
\texttt{A-n62-k8-99} & -35969 & 31.53 & 0.78 \\
\texttt{P-n65-k10-102} & -3 & 32.97 & 0.47 \\
\texttt{B-n68-k9-65} & -31001 & 34.95 & 1.28 \\
\texttt{P-n55-k15-88} & -2 & 35.49 & 0.18 \\
\texttt{A-n80-k10-14} & -105283 & 50.74 & 4.40 \\
\texttt{B-n66-k9-50} & -26520 & 55.13 & 0.51 \\
\texttt{G-n262-k25-316} & -1426535 & 58.47 & 3600$^\dagger$ \\
\texttt{E-n76-k14-102} & -1 & 87.94 & 0.68 \\
\texttt{A-n60-k9-57} & -1000 & 88.91 & 2.33 \\
\texttt{E-n76-k15-40} & -1 & 94.05 & 0.50 \\
\texttt{M-n121-k7-260} & -160097 & 179.5 & 3600$^\dagger$ \\
\texttt{M-n200-k17-12} & -121210 & 289.2 & 84.05 \\
\texttt{M-n151-k12-15} & -79996 & 331.0 & 102.6 \\
\texttt{E-n101-k14-158} & -3590 & 361.6 & 0.88 \\
\midrule
\multicolumn{4}{l}{\textbf{Roberti Set 3}} \\
\texttt{F-n45-k4\_a} & -14 & 0.28 & 0.05 \\
\texttt{P-n76-k4\_b} & -13 & 0.80 & 0.10 \\
\texttt{P-n76-k5\_b} & -14 & 0.89 & 0.08 \\
\texttt{E-n76-k8\_b} & -19 & 1.02 & 0.05 \\
\texttt{P-n101-k4\_b} & -13 & 1.25 & 0.58 \\
\texttt{F-n72-k4\_a} & 0 & 1.79 & 27.81 \\
\texttt{E-n101-k8\_b} & -19 & 2.44 & 0.17 \\
\texttt{E-n76-k7\_b} & -6 & 5.97 & 0.08 \\
\texttt{E-n76-k10\_b} & -8 & 7.58 & 0.05 \\
\texttt{P-n70-k10\_b} & -2 & 18.69 & 0.11 \\
\texttt{P-n76-k4\_a} & -3 & 18.98 & 0.61 \\
\texttt{P-n101-k4\_a} & -7 & 22.06 & 42.64 \\
\texttt{P-n76-k5\_a} & -4 & 27.82 & 1.78 \\
\texttt{E-n101-k8\_a} & -24 & 32.65 & 1.11 \\
\texttt{P-n70-k10\_a} & -3 & 35.98 & 0.07 \\
\texttt{E-n76-k7\_a} & -6 & 45.93 & 0.36 \\
\texttt{E-n76-k8\_a} & -7 & 52.88 & 0.26 \\
\texttt{E-n76-k10\_a} & -4 & 54.36 & 0.06 \\
\texttt{E-n76-k14\_a} & -4 & 54.81 & 0.04 \\
\texttt{E-n101-k14\_a} & -7 & 100.3 & 0.61 \\
\texttt{E-n76-k14\_b} & 0 & 101.2 & 0.61 \\
\texttt{E-n101-k14\_b} & 0 & 533.3 & 3.35 \\
\texttt{M-n151-k12\_b} & -4 & 1329.2 & 4.41 \\
\texttt{M-n121-k7\_a} & -15 & 2087.1 & 3600$^\dagger$ \\
\texttt{M-n121-k7\_b} & -9 & 3600$^\dagger$ & 42.54 \\
\texttt{M-n151-k12\_a} & -6 & 3600$^\dagger$ & 4.01 \\
\texttt{M-n200-k17\_b} & 0 & 3600$^\dagger$ & 649.1 \\
\texttt{M-n200-k17\_a} & -6 & 3600$^\dagger$ & 17.11 \\
\texttt{M-n200-k16\_a} & -5 & 3600$^\dagger$ & 19.47 \\
\texttt{M-n200-k16\_b} & -1 & 3600$^\dagger$ & 232.6 \\
\texttt{F-n135-k7\_a} & -14 & 3600$^\dagger$ & 14.77 \\
\end{longtable}
\endgroup

\begingroup\footnotesize
\setlength{\tabcolsep}{4pt}
\begin{longtable}{lrrrrrr}
\caption{Per-instance runtime (s) under each ablation configuration (300\,s limit; $\dagger$ = timeout). Columns match the ladder of Table~\ref{tbl:ablation}: GSEC only; +capacity cuts; +reduced-cost fixing; +comb/RGLM; +all-pairs bounds (+ap); +SPI.}\label{tbl:ablationinstances}\\
\toprule
Instance & GSEC & +cap & +fix & +comb & +ap & +SPI \\
\midrule
\endfirsthead
\multicolumn{7}{c}{\tablename\ \thetable\ -- continued from previous page}\\
\toprule
Instance & GSEC & +cap & +fix & +comb & +ap & +SPI \\
\midrule
\endhead
\midrule\multicolumn{7}{r}{\textit{continued on next page}}\\
\endfoot
\bottomrule
\endlastfoot
\midrule
\multicolumn{7}{l}{\textbf{SPPRCLIB}} \\
\texttt{A-n54-k7-149} & 300$^\dagger$ & 12.5 & 12.5 & 12.5 & 12.7 & 12.7 \\
\texttt{A-n60-k9-57} & 300$^\dagger$ & 88.3 & 88.0 & 79.3 & 88.0 & 88.0 \\
\texttt{A-n61-k9-80} & 23.5 & 11.1 & 11.3 & 11.4 & 11.3 & 11.3 \\
\texttt{A-n62-k8-99} & 300$^\dagger$ & 31.2 & 31.6 & 27.3 & 31.5 & 31.4 \\
\texttt{A-n63-k10-44} & 40.0 & 5.3 & 5.4 & 5.4 & 5.4 & 5.3 \\
\texttt{A-n63-k9-157} & 116.1 & 11.5 & 9.6 & 9.7 & 9.5 & 9.6 \\
\texttt{A-n64-k9-45} & 300$^\dagger$ & 26.0 & 26.4 & 25.9 & 26.4 & 26.2 \\
\texttt{A-n65-k9-10} & 19.6 & 6.3 & 6.4 & 6.3 & 6.3 & 6.3 \\
\texttt{A-n69-k9-42} & 31.9 & 4.0 & 4.1 & 4.0 & 4.0 & 4.0 \\
\texttt{A-n80-k10-14} & 63.8 & 50.2 & 50.4 & 50.8 & 50.6 & 50.9 \\
\texttt{B-n45-k6-54} & 6.4 & 5.5 & 5.4 & 5.5 & 5.5 & 5.5 \\
\texttt{B-n50-k8-40} & 300$^\dagger$ & 15.4 & 15.1 & 15.0 & 15.2 & 15.3 \\
\texttt{B-n52-k7-15} & 4.0 & 3.0 & 2.9 & 3.0 & 3.0 & 3.0 \\
\texttt{B-n57-k7-20} & 2.6 & 2.6 & 2.8 & 2.8 & 3.5 & 3.5 \\
\texttt{B-n66-k9-50} & 300$^\dagger$ & 54.9 & 54.3 & 46.8 & 55.0 & 54.5 \\
\texttt{B-n67-k10-26} & 201.0 & 2.8 & 2.8 & 2.8 & 2.9 & 2.9 \\
\texttt{B-n68-k9-65} & 300$^\dagger$ & 34.9 & 35.2 & 35.0 & 35.3 & 35.3 \\
\texttt{B-n78-k10-70} & 300$^\dagger$ & 23.0 & 23.3 & 23.3 & 23.5 & 23.3 \\
\texttt{E-n101-k14-158} & 300$^\dagger$ & 300$^\dagger$ & 300$^\dagger$ & 300$^\dagger$ & 300$^\dagger$ & 300$^\dagger$ \\
\texttt{E-n101-k8-291} & 50.0 & 25.0 & 28.3 & 28.3 & 28.7 & 28.4 \\
\texttt{E-n76-k10-72} & 32.1 & 12.8 & 13.0 & 13.1 & 12.9 & 12.9 \\
\texttt{E-n76-k14-102} & 189.1 & 105.4 & 88.1 & 85.4 & 99.2 & 98.6 \\
\texttt{E-n76-k15-40} & 114.2 & 66.1 & 94.7 & 95.0 & 96.6 & 95.8 \\
\texttt{E-n76-k7-44} & 21.6 & 6.7 & 6.9 & 6.0 & 7.0 & 6.9 \\
\texttt{G-n262-k25-316} & 139.4 & 56.9 & 58.5 & 57.3 & 56.9 & 56.2 \\
\texttt{M-n101-k10-97} & 9.4 & 18.7 & 18.7 & 18.7 & 18.8 & 18.5 \\
\texttt{M-n121-k7-260} & 201.8 & 176.5 & 179.5 & 178.5 & 184.2 & 177.0 \\
\texttt{M-n151-k12-15} & 300$^\dagger$ & 300$^\dagger$ & 300$^\dagger$ & 300$^\dagger$ & 300$^\dagger$ & 300$^\dagger$ \\
\texttt{M-n200-k16-143} & 11.3 & 12.2 & 11.8 & 11.9 & 79.2 & 80.3 \\
\texttt{M-n200-k17-12} & 158.2 & 294.9 & 288.7 & 289.2 & 300$^\dagger$ & 300$^\dagger$ \\
\texttt{P-n101-k4-174} & 30.8 & 16.5 & 16.7 & 16.6 & 17.7 & 18.0 \\
\texttt{P-n50-k10-24} & 10.6 & 7.4 & 10.7 & 10.8 & 10.7 & 10.8 \\
\texttt{P-n50-k7-92} & 17.8 & 14.0 & 13.9 & 12.5 & 13.8 & 13.9 \\
\texttt{P-n50-k8-19} & 2.7 & 1.4 & 1.4 & 1.4 & 1.5 & 1.5 \\
\texttt{P-n51-k10-30} & 25.6 & 20.8 & 20.9 & 20.9 & 20.6 & 20.8 \\
\texttt{P-n55-k10-44} & 29.1 & 21.4 & 20.7 & 20.5 & 20.4 & 20.6 \\
\texttt{P-n55-k15-88} & 24.0 & 34.2 & 35.2 & 35.3 & 35.1 & 34.8 \\
\texttt{P-n55-k7-116} & 11.1 & 2.3 & 2.3 & 2.3 & 2.3 & 2.3 \\
\texttt{P-n55-k8-260} & 23.9 & 8.9 & 11.3 & 11.3 & 11.1 & 11.4 \\
\texttt{P-n60-k10-24} & 10.0 & 2.2 & 2.3 & 2.2 & 2.3 & 2.3 \\
\texttt{P-n60-k15-8} & 32.9 & 24.1 & 23.4 & 23.2 & 24.0 & 24.4 \\
\texttt{P-n65-k10-102} & 64.0 & 37.1 & 32.7 & 36.9 & 32.8 & 32.6 \\
\texttt{P-n70-k10-12} & 2.4 & 1.1 & 1.1 & 1.2 & 1.2 & 1.2 \\
\texttt{P-n76-k4-41} & 9.2 & 8.3 & 8.3 & 8.2 & 8.6 & 8.6 \\
\texttt{P-n76-k5-16} & 1.6 & 1.3 & 1.4 & 1.3 & 1.9 & 1.9 \\
\midrule
\multicolumn{7}{l}{\textbf{Roberti Set 3}} \\
\texttt{E-n101-k14\_a} & 300$^\dagger$ & 101.2 & 101.1 & 101.1 & 101.0 & 100.4 \\
\texttt{E-n101-k14\_b} & 300$^\dagger$ & 300$^\dagger$ & 300$^\dagger$ & 300$^\dagger$ & 300$^\dagger$ & 300$^\dagger$ \\
\texttt{E-n101-k8\_a} & 34.0 & 32.1 & 32.8 & 32.3 & 32.7 & 33.3 \\
\texttt{E-n101-k8\_b} & 33.4 & 2.5 & 2.6 & 2.6 & 2.9 & 2.8 \\
\texttt{E-n76-k10\_a} & 300$^\dagger$ & 50.1 & 53.9 & 53.1 & 53.4 & 53.3 \\
\texttt{E-n76-k10\_b} & 55.5 & 7.8 & 7.7 & 7.7 & 7.6 & 7.6 \\
\texttt{E-n76-k14\_a} & 157.6 & 47.9 & 54.9 & 54.1 & 54.6 & 54.8 \\
\texttt{E-n76-k14\_b} & 211.7 & 117.9 & 101.0 & 102.2 & 119.0 & 117.4 \\
\texttt{E-n76-k7\_a} & 300$^\dagger$ & 46.0 & 46.3 & 69.8 & 45.8 & 45.3 \\
\texttt{E-n76-k7\_b} & 32.1 & 6.0 & 6.1 & 8.9 & 6.0 & 5.9 \\
\texttt{E-n76-k8\_a} & 300$^\dagger$ & 53.4 & 52.5 & 58.9 & 53.6 & 53.6 \\
\texttt{E-n76-k8\_b} & 5.1 & 1.0 & 1.0 & 1.0 & 1.1 & 1.1 \\
\texttt{F-n135-k7\_a} & 300$^\dagger$ & 300$^\dagger$ & 300$^\dagger$ & 300$^\dagger$ & 300$^\dagger$ & 300$^\dagger$ \\
\texttt{F-n45-k4\_a} & 1.0 & 0.3 & 0.3 & 0.3 & 0.4 & 0.4 \\
\texttt{F-n72-k4\_a} & 3.3 & 1.7 & 1.8 & 1.8 & 5.0 & 5.0 \\
\texttt{M-n121-k7\_a} & 300$^\dagger$ & 300$^\dagger$ & 300$^\dagger$ & 300$^\dagger$ & 300$^\dagger$ & 300$^\dagger$ \\
\texttt{M-n121-k7\_b} & 300$^\dagger$ & 300$^\dagger$ & 300$^\dagger$ & 300$^\dagger$ & 300$^\dagger$ & 300$^\dagger$ \\
\texttt{M-n151-k12\_a} & 300$^\dagger$ & 300$^\dagger$ & 300$^\dagger$ & 300$^\dagger$ & 300$^\dagger$ & 300$^\dagger$ \\
\texttt{M-n151-k12\_b} & 300$^\dagger$ & 300$^\dagger$ & 300$^\dagger$ & 300$^\dagger$ & 300$^\dagger$ & 300$^\dagger$ \\
\texttt{M-n200-k16\_a} & 300$^\dagger$ & 300$^\dagger$ & 300$^\dagger$ & 300$^\dagger$ & 300$^\dagger$ & 300$^\dagger$ \\
\texttt{M-n200-k16\_b} & 300$^\dagger$ & 300$^\dagger$ & 300$^\dagger$ & 300$^\dagger$ & 300$^\dagger$ & 300$^\dagger$ \\
\texttt{M-n200-k17\_a} & 300$^\dagger$ & 300$^\dagger$ & 300$^\dagger$ & 300$^\dagger$ & 300$^\dagger$ & 300$^\dagger$ \\
\texttt{M-n200-k17\_b} & 300$^\dagger$ & 300$^\dagger$ & 300$^\dagger$ & 300$^\dagger$ & 300$^\dagger$ & 300$^\dagger$ \\
\texttt{P-n101-k4\_a} & 86.8 & 21.7 & 21.5 & 22.0 & 23.0 & 22.6 \\
\texttt{P-n101-k4\_b} & 1.5 & 1.2 & 1.2 & 1.2 & 2.5 & 2.5 \\
\texttt{P-n70-k10\_a} & 107.0 & 34.8 & 36.0 & 35.9 & 35.7 & 35.1 \\
\texttt{P-n70-k10\_b} & 111.5 & 18.1 & 18.3 & 18.3 & 18.4 & 18.7 \\
\texttt{P-n76-k4\_a} & 61.0 & 18.8 & 19.0 & 19.3 & 19.0 & 19.3 \\
\texttt{P-n76-k4\_b} & 1.2 & 0.8 & 0.8 & 0.8 & 1.0 & 1.0 \\
\texttt{P-n76-k5\_a} & 104.1 & 27.7 & 27.4 & 28.3 & 27.7 & 27.6 \\
\texttt{P-n76-k5\_b} & 3.4 & 0.9 & 0.9 & 0.9 & 1.0 & 1.0 \\
\end{longtable}
\endgroup

\bibliographystyle{plainnat}
\bibliography{references}

\end{document}